\documentclass[11pt]{article}
\usepackage[utf8]{inputenc}

\usepackage{amsthm, amsfonts, amssymb, amsmath,enumitem, natbib, bbm, mathabx}
\usepackage{cases, amsthm}
\usepackage{fullpage}
\usepackage{mathtools}
\usepackage{mathrsfs} 
\usepackage{wasysym} 
\usepackage{verbatim}
\usepackage{fancyhdr}
\usepackage{graphicx}
\usepackage{bbm}
\usepackage{caption}
\usepackage{subcaption}
\usepackage{algorithm}
\usepackage{algorithmic}
\usepackage{authblk}
\usepackage{setspace}
\let\Algorithm\algorithm
\renewcommand\algorithm[1][]{\Algorithm[#1]\setstretch{1.6}}
\usepackage[colorlinks=true, pdfstartview=FitV,
linkcolor=blue, citecolor=blue, urlcolor=blue]{hyperref}
\usepackage{cleveref}

\newtheorem{definition}{Definition}[section]
\newtheorem{theorem}[definition]{Theorem}

\newtheorem{lemma}[definition]{Lemma}
\newtheorem{corollary}[definition]{Corollary}

\theoremstyle{remark}

\usepackage{accents}

\usepackage{xcolor}

\usepackage[top=1in, bottom=1in, left=1in, right=1in]{geometry}
\title{On the Missing Factor in Some Concentration
Inequalities for Martingales}

\author{Arun Kumar Kuchibhotla}
\affil{Department of Statistics \& Data Science, Carnegie Mellon University}
\date{}

\begin{document}
\maketitle
\begin{abstract}
    In this note, we improve some concentration inequalities for martingales with bounded increments. These results recover the missing factor in Freedman-style inequalities and are near optimal. We also provide minor refinements of concentration inequalities for functions of independent random variables. These proofs use techniques from the works of Bentkus and Pinelis.
\end{abstract}
\section{Introduction}\label{sec:introduction}
Concentration inequalities such as Hoeffding's inequality are central to modern (e.g., high-dimensional) statistical theory and inference. For sums of independent random variables that are mean zero and bounded, the classical Hoeffding, Bernstein, and Bennett inequalities have been improved in a series of works by Bentkus and Pinelis; see, for example,~\cite{bentkus2002remark,bentkus2004hoeffding,pinelis2006binomial,pinelis2014bennett}. These improvements recover the so-called ``missing factor'' in classical inequalities~\citep{talagrand1995missing}. Moreover, these inequalities are nearly optimal.

The goal of this article is to bring about similar improvements in martingale inequalities with bounded increments. In particular, we improve on Freedman's inequality and its refinements derived in~\cite{fan2012hoeffding,Fan2015}. To briefly mention the result, note that the traditional concentration inequalities are proved using $\mathbf{1}\{u \ge 0\} \le \exp(\lambda u)$ for any $\lambda \ge 0$. The improvements proposed in the works of Bentkus and Pinelis start with the inequality
\[
\mathbf{1}\{u \ge 0\} \le \left(1 + \frac{\lambda u}{\alpha}\right)_+^{\alpha} \le \exp(\lambda u)\quad\mbox{for all}\quad \lambda \ge 0, u\in\mathbb{R}.
\]
Hence, for any $x\in\mathbb{R}$ and any set of random variables $X_1, \ldots, X_n$,
\[
\mathbb{P}\left(\sum_{i=1}^n X_i \ge x\right) \le \inf_{\lambda \ge 0}\,\mathbb{E}\left[\left(1 + \frac{\lambda}{2}\left(\sum_{i=1}^n X_i - x\right)\right)_+^{\alpha}\right] \le \inf_{\lambda \ge 0}\,\mathbb{E}\left[\exp\left(\lambda\left(\sum_{i=1}^n X_i - x\right)\right)\right].
\]
The right hand side is the classical Cramer-Chernoff bound on the tail probability. The intermediate bound is the improvement suggested in~\cite{bentkus2002remark}, and we refer to these bounds as the Bentkus bound. There are two important differences between these bounds: (1) unlike the Cramer-Chernoff bound, the Bentkus bound is finite as long as the random variables have a finite $\alpha$-th moment; and (2) unlike the Cramer-Chernoff bound, Bentkus bound cannot be written as a product if $X_i$'s are mutually independent. This second aspect makes it difficult to work with the Bentkus bound in many cases. 

The results available for independent random variables are as follows.
\begin{enumerate}
    \item If $X_i\in [a_i, b_i], 1\le i\le n$ with $a_i \le 0 \le b_i$ and $\mathbb{E}[X_i] \le 0$, then for any $x\ge 0$,
    \[
    \mathbb{P}\left(\sum_{i=1}^n X_i \ge x\right) \le \inf_{\lambda \ge 0}\mathbb{E}\left[\left(1 + \lambda\left(\sum_{i=1}^n G_{1,i} - x\right)\right)_+\right] \le \inf_{\lambda \ge 0}\,\mathbb{E}\left[\exp\left(\lambda\left(\sum_{i=1}^n G_{1,i} - x\right)\right)\right],
    \]
    where $G_{1,i}\in\{a_i, b_i\}, 1\le i\le n$ are independent mean zero random variables. More formally, $G_{1,i}, 1\le i\le n$ are two-point random variables such that
    \begin{equation}\label{eq:Hoeffding-bounded}
    \mathbb{P}(G_{1,i} = a_i) = \frac{b_i}{b_i - a_i}\quad\mbox{and}\quad \mathbb{P}(G_{1,i} = b_i) = \frac{-a_i}{b_i - a_i}.
    \end{equation}
    \item If $X_i \le b_i, 1\le i\le n$ with $\mathbb{E}[X_i] \le 0$, $\mbox{Var}(X_i) \le \sigma_i^2$, then for any $x\ge0$,
    \[
    \mathbb{P}\left(\sum_{i=1}^n X_i \ge x\right) \le \inf_{\lambda \ge 0}\mathbb{E}\left[\left(1 + \frac{\lambda}{2}\left(\sum_{i=1}^n G_{2,i} - x\right)\right)_+^2\right] \le \inf_{\lambda \ge 0}\,\mathbb{E}\left[\exp\left(\lambda\left(\sum_{i=1}^n G_{2,i} - x\right)\right)\right],
    \]
    where $G_{2,i}\in\{-\sigma_i^2/b_i, b_i\}, 1\le i\le n$ are independent mean zero random variables. More formally, $G_{2,i}, 1\le i\le n$ are two-point random variables such that
    \[
    \mathbb{P}(G_{2,i} = -\sigma_i^2/b_i) = \frac{b_i^2}{\sigma_i^2 + b_i^2}\quad\mbox{and}\quad \mathbb{P}(G_{2,i} = b_i) = \frac{\sigma_i^2}{\sigma_i^2 + b_i^2}.
    \]
\end{enumerate}
Note that in both cases, we have lower bounds in terms of $G_{1,i}$'s and $G_{2,i}$'s, i.e., 
\[
\sup_{X_i\in[a_i,b_i], \mathbb{E}[X_i] = 0}\mathbb{P}\left(\sum_{i=1}^n X_i \ge x\right) \ge \mathbb{P}\left(\sum_{i=1}^n G_{1,i} \ge x\right),
\]
and
\[
\sup_{X_i\le b_i, \mathbb{E}[X_i] = 0, \mbox{Var}(X_i) \le \sigma_i^2}\mathbb{P}\left(\sum_{i=1}^n X_i \ge x\right) \ge \mathbb{P}\left(\sum_{i=1}^n G_{2,i} \ge x\right).
\]
The near-optimality of Bentkus bounds is that for $\alpha\in\{1, 2\}$, and $x$ in the support of $\sum_{i=1}^n G_{\alpha,i}$, we have
\begin{equation}\label{eq:optimality}
\inf_{\lambda \ge 0}\mathbb{E}\left[\left(1 + \frac{\lambda}{\alpha}\left(\sum_{i=1}^n G_{\alpha,i} - x\right)\right)_+^\alpha\right] ~\le~ \frac{e^{\alpha}}{\alpha}\mathbb{P}\left(\sum_{i=1}^n G_{\alpha,i} \ge x\right).
\end{equation}
Therefore, the Bentkus bounds presented above cannot be improved by more than a factor of $e^{\alpha}/\alpha$. It should be noted that the left-hand side of~\eqref{eq:optimality} is at most one for any $x \ge 0$ (e.g. by taking $\lambda = 0$), while the right hand side can be larger than $1$. Moreover, we note that the optimal Cramer-Chernoff bounds presented above are known in~\cite{bennett1962probability} and~\cite{hoeffding1963probability}. For readers unfamiliar with the Bentkus bounds,~\cite{bentkus2002remark} provides an accessible introduction to this literature. As noted in~\cite{hoeffding1963probability} and~\cite{bentkus2004hoeffding}, the optimal Cramer-Chernoff bound is better than the classical Bernstein, Bennett, and Prokhorov inequalities. The Bentkus bound is also better than Talagrand's inequality~\citep{talagrand1995missing}.

One more point to note in comparison with classical concentration inequalities is that the Bentkus bound is not just a function of the variance of the sum, but is a nontrivial function of individual variances. Using some domination inequalities for sums of Bernoulli random variables~\citep[Lemma 4.5]{bentkus2004hoeffding}, one can obtain (less precise) Bentkus bounds that depend only on $\sum_{i=1}^n (b_i - a_i)^2$ and $\sum_{i=1}^n \sigma_i^2$, respectively. For example,~\cite{Pinelis2006} proves
\[
\sup_{X_i\in[a_i, b_i], \mathbb{E}[X_i] \le 0}\,\mathbb{P}\left(\sum_{i=1}^n X_i \ge x\right) \le \inf_{\lambda \ge 0}\,\mathbb{E}\left[\left(1 + \frac{\lambda}{5}(v Z - x)\right)_+^5\right],\quad\mbox{where}\quad v^2 = \sum_{i=1}^n (b_i - a_i)^2,
\]
and~\cite{bentkus2004hoeffding} proves
\[
\sup_{X_i \le b_i, \mathbb{E}[X_i] \le 0, \mbox{Var}(X_i) \le \sigma_i^2}\,\mathbb{P}\left(\sum_{i=1}^n X_i \ge x\right) \le \inf_{\lambda \ge 0}\,\mathbb{E}\left[\left(1 + \frac{\lambda}{2}\left(\sum_{i=1}^n G_i' - x\right)\right)_+^2\right],
\]
where $G_i'\in\{-v^2/b, b\}$ with $b = \max_{1\le i\le n}b_i$, and $v^2 = \sum_{i=1}^n \sigma_i^2$. These inequalities have been extended to super-martingales in the works cited above, but with the restrictive assumption that $\sigma_i^2$ represents an almost sure upper bound on the conditional variance of $X_i$ given $\mathcal{F}_{i-1}$. This makes it difficult to obtain concentration inequalities for functions of independent random variables via the Doob martingale. The goal of this note is to obtain results without assuming that the conditional variances are almost surely bounded. 

\paragraph{On the Proof Techniques.} Most of the martingale analogues of the classical concentration inequalities (for independent data) are obtained using the Cramer-Chernoff technique. Some of the inequalities obtained time uniformly (i.e., simultaneously valid for all sample sizes) have also been obtained using the Cramer-Chernoff technique. The main idea in these proofs is the construction of an exponential supermartingale and an application of Ville's or Doob's maximal inequality~\citep{howard2021time}. The construction of this exponential supermartingale is based on bounds on the moment generating function. Unfortunately, there is no such simple supermartingale in relation to the Bentkus style bounds, and it is not obvious if there is a general way to extend/improve all the existing Cramer-Chernoff results to the Bentkus style bounds. In fact, it is not clear whether such analogues hold true, but we know that if such an analogue exists, then it is a significant improvement over the Cramer-Chernoff ones.  

Our proofs are closely related to the Lindeberg swapping technique commonly used in proving central limit theorems~\citep{lalley2014martingale}. In the literature on concentration inequalities, the closest result is Lemma 5.2 of~\cite{pinelis1998optimal}. Our results arise as direct corollaries of the proof of this lemma. This note exists only because these results have not appeared in the literature in this form before and might be of interest to the broader community.
\section{Improving Freedman-style Inequality}\label{sec:Freedman-inequality}
Suppose $(X_i, \mathcal{F}_i)_{i\ge0}$ with $X_0 = 0$ forms a (super)martingale difference sequence bounded above by $1$, i.e., $X_i$ is $\mathcal{F}_i$-measurable, and $\mathbb{E}[X_i|\mathcal{F}_{i-1}] \le 0$. Let $\sigma_i^2 := \mathbb{E}[X_i^2|\mathcal{F}_{i-1}]$. Consider the event
\[
\mathcal{E}_k := \left\{\sum_{i=1}^k X_i \ge x\quad\mbox{and}\quad \sum_{i=1}^k \sigma_i^2 \le v^2\right\}.
\]
For any $x \ge 0$ and $0 < v < \infty$, \cite{freedman1975tail} proved inequalities for $\mathbb{P}(\cup_{k\ge1}\mathcal{E}_k)$ that can be considered as extensions of the classical Bernstein inequality for independent random variables. \cite{Fan2015} improved Freedman inequalities by proving
\begin{equation}\label{eq:Fan-inequality}
    \mathbb{P}\left(\bigcup_{k=1}^n \mathcal{E}_k\right) \le \inf_{\lambda \ge 0}\,\mathbb{E}\left[\exp\left(\lambda\left(\sum_{i=1}^n G_i - x\right)\right)\right],
\end{equation}
where $G_i\in\{-v^2/n, 1\}, 1\le i\le n$ are mutually independent random variables with mean zero, i.e., 
\begin{equation}\label{eq:optimal-Freedman}
\mathbb{P}(G_i = -v^2/n) = \frac{1}{1 + v^2/n}\quad\mbox{and}\quad \mathbb{P}(G_i = 1) = \frac{v^2/n}{1 + v^2/n}.
\end{equation}
(Although Theorem 2.1 of~\cite{Fan2015} is not presented in this form, it is equivalent to the form above; see Eq. (41) of~\cite{Fan2015}.) Moreover,~\cite{Fan2015} show that this bound improves many results as described in Remark 2.1 and Corollary 2.2 there. 

The following is the main result of this section. 
\begin{theorem}\label{thm:Freedman-improvement}
    Suppose $(X_i, \mathcal{F}_i), i\ge 0$ forms a supermartingale difference sequence with $X_0 = 0$ and $\mathbb{P}(X_i \le 1) = 1$\footnote{We assume the bound to be $1$, without loss of generality.} Then
    \begin{equation}\label{eq:finite-n-Freedman}
    \mathbb{P}\left(\bigcup_{k=1}^n \mathcal{E}_k\right) \le \inf_{\lambda \ge 0}\,\mathbb{E}\left[\left(1 + \frac{\lambda}{2}\left(\sum_{i=1}^n G_i - x\right)\right)_+^2\right] \le \inf_{\lambda \ge 0}\,\mathbb{E}\left[\left(1 + \frac{\lambda}{2}(\widetilde{\Pi}_{v^2} - x)\right)_+^2\right],
    \end{equation}
    where $G_i$'s are independent and identically distributed random variables with distribution specified in~\eqref{eq:optimal-Freedman}, and $\widetilde{\Pi}_{v^2} \sim \mathrm{Poisson}(v^2) - v^2$. In particular, we obtain
    \begin{equation}\label{eq:infinite-Freedman}
    \begin{split}
    \mathbb{P}\left(\sum_{i=1}^n X_i \ge x\quad\mbox{and}\quad \sum_{i=1}^n \sigma_i^2 \le v^2\quad\mbox{for some}\quad n\ge 1\right) &\le \inf_{\lambda \ge 0}\,\mathbb{E}\left[\left(1 + \frac{\lambda}{2}(\widetilde{\Pi}_{v^2} - x)\right)_+^2\right]\\
    &\le \frac{e^2}{2}\mathbb{P}^{\circ}(\widetilde{\Pi}_{v^2} \ge x),
    \end{split}
    \end{equation}
    where $x\mapsto\mathbb{P}^{\circ}(\widetilde{\Pi}_{v^2} \ge x)$ is the least log-concave majorant of $x\mapsto\mathbb{P}(\widetilde{\Pi}_{v^2} \ge x).$
\end{theorem}
To our knowledge, Theorem~\ref{thm:Freedman-improvement} is new. The second inequality of~\eqref{eq:finite-n-Freedman} follows from Proposition 2.8 of~\cite{pinelis2014bennett}. Inequality~\eqref{eq:infinite-Freedman} follows from the second inequality of~\eqref{eq:finite-n-Freedman} and Lemma 4.2 of~\cite{bentkus2004hoeffding}. The first inequality of~\eqref{eq:finite-n-Freedman} is already near optimal for independent data because it is exactly the Bentkus bound for sums of IID mean zero random variables that are bounded above. Moreover, from the works cited in Section~\ref{sec:introduction}, the first inequality of~\eqref{eq:finite-n-Freedman} is an improvement of Theorem 2.1 of~\cite{Fan2015}, and therefore an improvement of the Freedman and de la Pena inequalities. In particular, Theorem~\ref{thm:Freedman-improvement} recovers the missing factor similar to how Bentkus' inequality for independent data recovers the missing factor in Bernstein's inequality.

The proof of Theorem~\ref{thm:Freedman-improvement} is provided in Section~\ref{appsec:Freedman-improvement}. Traditional proofs of Freedman's inequality and its improvements are based on the exponential function, and the fact that the exponential of a sum is the product of exponentials significantly helps the analysis. Our proof is a repeated application of Lemma 4.4(i) and Lemma 4.5 of~\cite{bentkus2004hoeffding}. 

\section{Improving Azuma-Hoeffding Inequality}\label{sec:Azuma-Hoeffding}
In Section~\ref{sec:Azuma-Hoeffding}, we considered the case of supremartingale differences that are bounded above by a non-random constant. In this section, we obtain an inequality that allows for a random upper bound and yields a generalization of the Azuma-Hoeffding inequality. Suppose $(X_i, \mathcal{F}_i)_{i\ge0}$ is a supermartingale difference sequence such that for some $\mathcal{F}_{i-1}$-measurable random variables $B_i$, we have $\mathbb{P}(X_i \le B_i|\mathcal{F}_{i-1}) = 1$ and $\sigma_i^2 = \mathbb{E}[X_i^2|\mathcal{F}_{i-1}]$. Define
\[
s_i := \frac{1}{2}\left(B_i + \frac{\sigma_i^2}{B_i}\right) \quad\mbox{and}\quad V_k^2 = \sum_{i=1}^k s_i^2.
\]
Consider the event
\[
\mathcal{E}_k = \left\{\sum_{i=1}^k X_i \ge x\;\mbox{and}\; V_k^2 \le v^2\right\}.
\]
With this notation, we have the following result.
\begin{theorem}\label{thm:Azuma-Hoeffding}
    For any $n\ge1$, $x\ge 0$, and $0 < v < \infty$, with $Z\sim N(0, 1)$,
    \[
    \mathbb{P}\left(\bigcup_{k=1}^n \mathcal{E}_k\right) \le \inf_{\lambda \ge 0}\,\mathbb{E}\left[\left(1 + \frac{\lambda}{5}(vZ - x)\right)_+^5\right] \le 5!(e/5)^5\mathbb{P}(vZ \ge x).
    \]
    In particular, 
    \[
    \mathbb{P}\left(\sum_{i=1}^n X_i \ge x\quad\mbox{and}\quad \sum_{i=1}^n s_i^2 \le v^2\mbox{ for some }n\ge1\right) \le \inf_{\lambda \ge 0}\,\mathbb{E}\left[\left(1 + \frac{\lambda}{5}(vZ - x)\right)_+^5\right].
    \]
\end{theorem}
Theorem~\ref{thm:Azuma-Hoeffding} improves inequality (2.22) of~\cite{Fan2015} and fills the gap in the literature; see the discussion following Eq. (2.22) of~\cite{Fan2015}. The proof of Theorem~\ref{thm:Azuma-Hoeffding} follows a structure similar to that of Theorem~\ref{thm:Freedman-improvement} using Lemma 5.1.3 of~\cite{Pinelis2006}.

Note that if $\mathbb{P}(A_i \le X_i \le B_i|\mathcal{F}_{i-1}) = 1$ for some non-negative $\mathcal{F}_{i-1}$-measurable random variables $A_{i}, B_i$, then $\sigma_i^2 \le -A_iB_i$ and hence $s_i \le (B_i - A_i)/2$. This fact combined with Theorem~\ref{thm:Azuma-Hoeffding} implies the following corollary.
\begin{corollary}\label{cor:Azuma-Hoeffding}
    Suppose $(X_i, \mathcal{F}_i)_{i\ge0}$ is a supermartingale difference sequence such that $\mathbb{P}(A_i \le X_i \le B_i|\mathcal{F}_{i-1}) = 1$ almost surely. Then
    \[
    \mathbb{P}\left(\sum_{i=1}^n X_i \ge x\quad\mbox{and}\quad \sum_{i=1}^n (B_i - A_i)^2 \le 4v^2\mbox{ for some }n\ge1\right) \le \inf_{\lambda \ge 0}\,\mathbb{E}\left[\left(1 + \frac{\lambda}{5}(vZ - x)\right)_+^5\right].
    \]
\end{corollary}
Corollary~\ref{cor:Azuma-Hoeffding} improves the classical Azuma-Hoeffding and McDiarmid inequalities. It also improves Theorems 2.5 and 2.6 of~\cite{van2002hoeffding}.
\section{Unbounded Martingale Increments and Functions of Independent Random Variables}
The results in previous sections are concerned with bounded supermartingale differences, and in many practical settings, one often encounters differences that are unbounded but bounded with high probability, for example, via certain tail assumption on the differences. This leads to generalizations of McDiarmid's inequality with unbounded influences. There is a thriving literature on this problem~\cite{li2023distribution,li2024concentration,marchina2018concentration,maurer2021concentration,kutin2002extensions,maurer2019bernstein,kontorovich2014concentration}. As with the classical concentration inequalities, most of these are based on the moment generating function and Cramer-Chernoff technique. In the following, we provide some results of the Bentkus type for these problems. Some results in this direction are provided in~\cite{marchina2018concentration}.

The following results are based on the framework and closely related to Theorem 2.1 of~\cite{marchina2018concentration}, which itself is based on the results of~\cite{bentkus2008extension,bentkus2010bounds}. To present the main result, we need some notation and preliminary results. For any two random variables $U$ and $V$,  we say $U$ is stochastically dominated by $V$ of order $\alpha \ge 0$ (denoted by $\preceq_\alpha$) if
\[
\mathbb{E}[(U - t)_{+}^{\alpha}] \le \mathbb{E}[(V - t)_+^{\alpha}]\quad\mbox{for all}\quad t\in\mathbb{R}.
\]
In particular, $U \preceq_{0} V$ is equivalent to $\mathbb{P}(U > t) \le \mathbb{P}(V > t)$ for all $t\in\mathbb{R}$. Note that $U \preceq_{\alpha} V$ is not equivalent to $-V \preceq_{\alpha} -U$ in general, unless $V$ is absolutely continuous with respect to the Lebesgue measure.
\begin{lemma}\label{lem:Hoeffding-unbounded}
Suppose $X\preceq_1 W$ and $-X \preceq_1 -T$ for some random variables $T$ and $W$ satisfying $T \preceq_0 W$. Then $\mathbb{E}[T] \le \mathbb{E}[X] \le \mathbb{E}[W]$. Moreover, for any $q\in[0, 1]$, if $a_q := \inf\{x:\, \mathbb{P}(T \le x) \ge 1 - q\}$ and $b_q := \inf\{x:\,\mathbb{P}(W \le x) \ge 1 - q\}$, then the random variable $\xi_q$ defined by
\[
\mathbb{P}(\xi_q \le x) = \begin{cases}\mathbb{P}(T \le x), &\mbox{if }x < a_q,\\
1-q, &\mbox{if }a_q \le x < b_q,\\
\mathbb{P}(W \le x), &\mbox{if }x \ge b_q,\end{cases}
\]
satisfies the following:
\begin{enumerate}
    \item $\xi_0 \overset{d}{=} T$, $\xi_1 \overset{d}{=} W$, $\mathbb{E}[T] \le \mathbb{E}[\xi_q] \le \mathbb{E}[W]$ for all $q\in[0, 1]$.
    \item $q\mapsto \mathbb{E}[\xi_q]$ is a differentiable non-decreasing function.
    \item If $q_0$ is the largest number such that $\mathbb{E}[\xi_{q_0}] = \mathbb{E}[X]$, then $X\preceq_1\xi_{q_0}$, i.e., $\mathbb{E}[(X - t)_+] \le \mathbb{E}[(\xi_{q_0} - t)_+]$ for all $t\in\mathbb{R}$.
\end{enumerate}
\end{lemma}
This lemma is an expanded version of Lemma 4.3 of~\cite{marchina2018concentration} and follows from Theorem 1 of~\cite{bentkus2010bounds}. In what follows, we are concerned with mean zero random variables and, therefore, we propose the notation $\xi_{T, W}$ for $\xi_{q_0}$ in Lemma~\ref{lem:Hoeffding-unbounded} with $q_0$ chosen so that $\mathbb{E}[\xi_{q_0}] = 0$. The reader should note that $\xi_q$ depends on $T$ and $W$ only through their distributions. 

The following result is one of the main results in this section. We remark/stress here that this is a known result and essentially follows by combining Theorem 1 of~\cite{bentkus2010bounds} and the basic proof techniques of~\cite{bentkus2008extension,bentkus2004hoeffding}. Given the importance of such a result for applications in statistics and machine learning, this result needs to be spelled out explicitly. 
\begin{theorem}\label{thm:martingale-stoch-bound}
    Suppose $(X_i,\mathcal{F}_{i})_{i\ge0}$ is a real-valued martingale difference sequence and suppose there exist random variables $T_i\preceq_0 W_i$ such that $X_i \preceq_1 W_i$ and $-X_i \preceq -T_i$ (conditional on $\mathcal{F}_{i-1}$), i.e., with probability 1,
    \[
    \mathbb{E}[(X_i - t)_+|\mathcal{F}_{i-1}] \le \mathbb{E}[(T_i - t)_+]\quad\mbox{and}\quad \mathbb{E}[(t - X_i)_+|\mathcal{F}_{i-1}] \le \mathbb{E}[(t - W_i)_+]\quad\mbox{for all}\quad t\in\mathbb{R}. 
    \]
    Then, for mutually independent random variables $\xi_{T_i, W_i}, 1\le i\le n$, we have
    \begin{equation}\label{eq:positive-1}
    \mathbb{E}\left[\left(\sum_{i=1}^n X_i - t\right)_+\right] \le \mathbb{E}\left[\left(\sum_{i=1}^n \xi_{T_i, W_i} - t\right)_+\right]\quad\mbox{for all}\quad t\in\mathbb{R}.
    \end{equation}
    In addition, if $T_i, W_i$ are absolutely continuous with respect to the Lebesgue measure, then 
    \begin{equation}\label{eq:negative-1}
    \mathbb{E}\left[\left(t-\sum_{i=1}^n X_i\right)_+\right] \le \mathbb{E}\left[\left(t + \sum_{i=1}^n \xi_{-W_i, -T_i}\right)_+\right]\quad\mbox{for all}\quad t\in\mathbb{R}.
    \end{equation}
    Consequently, for any $x\ge 0$,
    \begin{equation}\label{eq:martingale-inequalities}
        \begin{split}
            \mathbb{P}\left(\sum_{i=1}^n X_i \ge x\right) &\le \inf_{\lambda \ge 0}\,\mathbb{E}\left[\left(1 + \lambda\left(\sum_{i=1}^n \xi_{T_i, W_i} - x\right)\right)_+\right],\\
            \mathbb{P}\left(\sum_{i=1}^n X_i \le -x\right) &\le \inf_{\lambda \ge 0}\,\mathbb{E}\left[\left(1 + \lambda\left(\sum_{i=1}^n \xi_{-W_i, -T_i} - x\right)\right)_+\right],
        \end{split}
    \end{equation}
\end{theorem}
It should be stressed here that $\xi_{T_i, W_i}, 1\le i\le n$ are independent and $T_i, W_i$ are independent of $\mathcal{F}_n.$ The proof of Theorem~\ref{thm:martingale-stoch-bound} is omitted because it easily follows from the proofs in~\cite{bentkus2008extension}. Moreover, inequalities~\eqref{eq:positive-1} and~\eqref{eq:negative-1} can be combined to yield an improved tail bound in~\eqref{eq:martingale-inequalities} by noting that if $\mathbb{E}[(R_1 - t)_+] \le \mathbb{E}[(R_2 - t)_+]$ for all $t\in\mathbb{R}$ and $\mathbb{E}[R_1] = \mathbb{E}[R_2]$, then $\mathbb{E}[(t - R_1)_+] \le \mathbb{E}[(t - R_2)_+]$ for all $t\in\mathbb{R}$; see~\citet[Proposition 3]{bentkus2008extension}.
This result yields the following important corollary for functions of independent random variables.
\begin{corollary}\label{cor:functions-of-ind}
    Suppose $X_1, \ldots, X_n$ are independent random variables taking values in some space $\mathcal{X}$ and $Z = F(X_1, \ldots, X_n)$ for some function $F:\mathcal{X}^n\to\mathbb{R}$. For $1\le i\le n$, let 
    \[
    Z_i = F(X_1, \ldots, X_{i-1}, X_i', X_{i+1}, \ldots, X_n),\quad\mbox{and}\quad \Delta_i := Z - Z_i,
    \]
    where $X_i'$ is an independent copy of $X_i$ (independent of $X_1,\ldots, X_n$). Suppose that we have random variables $T_i\preceq_0 W_i$ such that
    \begin{equation}\label{eq:interaction-stoc-bnd}
    \mathbb{E}[(\Delta_i - t)_+|\mathcal{F}_{i-1}] \le \mathbb{E}[(W_i - t)_+]\quad\mbox{and}\quad \mathbb{E}[(t - \Delta_i)_+|\mathcal{F}_{i-1}] \le \mathbb{E}[(t - T_i)_+].
    \end{equation}
    Then
    \[
    \mathbb{E}[(Z - \mathbb{E}[Z] - t)_+] \le \mathbb{E}\left[\left(\sum_{i=1}^n \xi_{T_i, W_i} - t\right)_+\right]\quad\mbox{for all}\quad t\in\mathbb{R}.
    \]
    Consequently, for all $x\in\mathbb{R}$,
    \[
    \mathbb{P}\left(Z - \mathbb{E}[Z] \ge x\right) \le \inf_{\lambda\ge0}\,\mathbb{E}\left[\left(1 + \lambda\left(\sum_{i=1}^n \xi_{T_i, W_i} - x\right)\right)_+\right].
    \]
\end{corollary}
Corollary~\ref{cor:functions-of-ind} follows from Theorem~\ref{thm:martingale-stoch-bound} by taking $\mathcal{F}_i = \sigma(\{X_1, \ldots, X_i\})$ and $X_i = \mathbb{E}[\Delta_i|\mathcal{F}_{i}]$ which is mean zero and forms a martingale difference sequence with respect to $\mathcal{F}_i$. More formally, the Doob martingale decomposition yields
\[
Z - \mathbb{E}[Z] ~=~ \sum_{i=1}^n (\mathbb{E}[Z|\mathcal{F}_i] - \mathbb{E}[Z|\mathcal{F}_{i-1}]) ~=~ \sum_{i=1}^n \mathbb{E}[\Delta_i|\mathcal{F}_i],
\]
where $\mathcal{F}_0$ is the trivial $\sigma$-field. This shows that condition~\eqref{eq:interaction-stoc-bnd} can be relaxed by replacing $\Delta_i$ with $\mathbb{E}[\Delta_i|\mathcal{F}_i]$. Given the use of the (first-order) influence in traditional concentration inequalities for functions of independent random variables, we feel that the condition~\eqref{eq:interaction-stoc-bnd} is more natural and intuitive. 
In Theorem~\ref{thm:martingale-stoch-bound} and Corollary~\ref{cor:functions-of-ind}, if $T_i$ and $W_i$ are taken to be degenerate random variables, then they yield refinements of the Azuma-Hoeffding and McDiarmid inequalities, respectively. This holds because for the degenerate $T_i, W_i$, $\xi_{T_i, W_i}$ is the two-point random variable presented in~\eqref{eq:Hoeffding-bounded}. In the special case where $-T_i \overset{d}{=} W_i \overset{d}{=} B_i$ for some constant $B_i > 0$, we get $\xi_{T_i, W_i} \overset{d}{=} B_i\varepsilon_i$ for a Rademacher random variable $\varepsilon_i$ independent of $B_i$ (i.e., $\mathbb{P}(\varepsilon_i = -1) = \mathbb{P}(\varepsilon_i=1) = 1/2$). Corollary~\ref{cor:functions-of-ind} should be compared with Theorems 3 \& 4 of~\cite{maurer2021concentration}, Lemma 3 \& Theorem 4 of~\cite{li2024concentration}, and Theorems 3.4 \& 3.5 of~\cite{li2023distribution}. Unlike the results of these papers, our result is derived without making any specific assumptions about the tail behavior of $T_i, W_i$. For example, one can allow $T_i, W_i$ to have tails that decay as slowly as quadratically. On the other hand, we require $\Delta_i$ to be stochastically bounded conditional on $\mathcal{F}_{i-1}$ for each $i$ individually, while in contrast,~\cite{maurer2021concentration} only require the average (over $i$) conditional sub-Gaussianity/sub-exponentially constant to be bounded. 

In the following, we discuss a few applications of Corollary~\ref{cor:functions-of-ind} to show its generality.
\paragraph{Sums of Random Variables.} Suppose $X_1, \ldots, X_n$ are independent random variables taking values in a Banach space $(\mathbb{B}, \|\cdot\|)$ for some norm $\|\cdot\|$. Take
\[
Z = \left\|\sum_{i=1}^n X_i\right\|.
\]
Note that we do not require $X_i$'s to be centered. By the triangle inequality, we have
\[
|\Delta_i| \le \|X_i - X_i'\|\quad\mbox{for all}\quad 1\le i\le n.
\]
Hence, we can take $-T_i = W_i = \|X_i - X_i'\|$. Corollary~\ref{cor:functions-of-ind} implies that for all $x\in\mathbb{R}$
\begin{equation}\label{eq:norms-of-sums}
\mathbb{P}\left(\left\|\sum_{i=1}^n X_i\right\| - \mathbb{E}\left\|\sum_{i=1}^n X_i\right\| \ge x\right) \le \inf_{\lambda \ge 0}\mathbb{E}\left[\left(1 + \lambda\left(\sum_{i=1}^n \varepsilon_i\|X_i - X_i'\| - x\right)\right)_+\right]\quad\mbox{for all}\quad x\in\mathbb{R}.
\end{equation}
If one has a tail assumption such as $\mathbb{P}(\|X_i - X_i'\| \ge t) \le \overline{F}_i(t)$ for some non-increasing function $\overline{F}(\cdot)$, then the right hand side can be further bounded by 
\[
\inf_{\lambda \ge 0}\mathbb{E}\left[\left(1 + \lambda\left(\sum_{i=1}^n \varepsilon_iW_i - x\right)\right)_+\right],
\]
where $W_i, 1\le i\le n$ are independent random variables whose survival functions are $\overline{F}_i(t)$. For instance, if $\|X_i - X_i'\|$ is $\sigma_i$-subGaussian, i.e., $\mathbb{P}(\|X_i - X_i'\| \ge t) \le \exp(-t^2/(2\sigma_i^2))$ for all $t\ge 0$ and some $\sigma_i > 0$, then we can take $W_i$'s to be Weibull random variables with shape parameter $2$. 

Inequality~\eqref{eq:norms-of-sums} should be compared with Proposition 7 of~\cite{maurer2021concentration} and Theorem 35 of~\cite{li2024concentration}.
\paragraph{Lipschitz Functions.} Suppose $X_1, \ldots, X_n$ are independent random variables taking values in a Banach space $(\mathbb{B}, \|\cdot\|)$ for some norm $\|\cdot\|$. Take $Z = F(X_1, \ldots, X_n)$ for some function $F$ satisfying
\[
|F(x_1, \ldots, x_n) - F(y_1, \ldots, y_n)| \le \sum_{i=1}^n L_i\|x_i - y_i\|\quad\mbox{for all}\quad x_i, y_i.
\]
Then $|\Delta_i| \le L_i\|X_i - X_i'\|$ for all $1\le i\le n$ and hence, for all $x\in\mathbb{R}$,
\[
\mathbb{P}\left(\left\|\sum_{i=1}^n X_i\right\| - \mathbb{E}\left\|\sum_{i=1}^n X_i\right\| \ge x\right) \le \inf_{\lambda \ge 0}\mathbb{E}\left[\left(1 + \lambda\left(\sum_{i=1}^n L_i\varepsilon_i\|X_i - X_i'\| - x\right)\right)_+\right].
\]
This inequality should be compared with Theorem 5 of~\cite{li2024concentration}.


\section{Comments and Conjectures for Unbounded Martingale Increments}\label{sec:conjectures}
Theorem~\ref{thm:martingale-stoch-bound} provides tail bounds for martingales with unbounded martingale increments, but requires stochastic dominance of the increments conditional on their past. This can be a restrictive assumption, and additionally, the bounds do not account for variances. One way to tackle both of these issues is to consider winsorized random variables, similar to Corollary 2.4 of~\cite{pinelis2006binomial}. 

For example, if $(X_i, \mathcal{F}_i)_{i\ge0}$ is a real-valued martingale difference sequence, then for any $y > 0$, $(\min\{X_i, y\}, \mathcal{F}_i)_{i\ge0}$ is a super-martingale difference sequence and hence, setting $\sigma_i^2 = \mathbb{E}[X_i^2|\mathcal{F}_{i-1}]$, Theorem~\ref{thm:Freedman-improvement} implies
\begin{equation}\label{eq:truncation-and-Thm2.1}
\begin{split}
&\mathbb{P}\left(\sum_{i=1}^k X_i \ge x\quad\mbox{and}\quad\sum_{i=1}^k \sigma_i^2 \le v^2\mbox{ for some }1\le k\le n\right)\\ 
&\le \inf_{\lambda \ge 0}\,\mathbb{E}\left[\left(1 + \frac{\lambda}{2}(y\widetilde{\Pi}_{v^2/y} - x)\right)_+^2\right] + \mathbb{P}\left(\max_{1\le k\le n}X_i > y\right),\quad\mbox{for any}\quad y \ge 0.
\end{split}
\end{equation}
We can replace $\sigma_i^2$ with $\sigma_i^2(y) = \mathbb{E}[(\min\{X_i, y\})^2|\mathcal{F}_{i-1}]$. 
From Theorem~\ref{thm:Freedman-improvement}, it is clear that one can replace the centered Poisson, $\widetilde{\Pi}_{v^2/y}$, in~\eqref{eq:truncation-and-Thm2.1} with a sum of $n$ independent identically distributed two-point random variables. The main result of~\cite{hahn1997approximation} suggests that if $X_i$'s are not too heterogeneously distributed, then for the optimal choice of $y$, inequality~\eqref{eq:truncation-and-Thm2.1} may be unimprovable except for a constant multiplier. Finally, inequality~\eqref{eq:truncation-and-Thm2.1} can be considered as an improvement of Corollary 2.2 of~\cite{fan2012hoeffding}.

Allowing the winsorization to be different for the martingale increments, one can get the following result. Following Theorem 2.3 of~\cite{pinelis2014optimal}, for any random variable $U$, define
\[
Q_{\alpha}(U; \delta) = \inf_{t\in\mathbb{R}}\, t + \frac{(\mathbb{E}[(U - t)_+^{\alpha}])^{1/\alpha}}{\delta^{1/\alpha}},\quad\mbox{for}\quad \alpha\in(0,\infty), \delta\in(0, 1).
\]
Note that $Q_{\alpha}(rU; \delta) = rQ_{\alpha}(U; \delta)$ for any $r \ge 0$.
\begin{corollary}\label{cor:different-truncations}
    Suppose $(X_i, \mathcal{F}_i)_{i\ge0}$ is a martingale difference sequence such that, for non-increasing, non-negative function $t\mapsto g(t)$,
    \begin{equation}\label{eq:survival-assumption}
    \mathbb{E}\left[\left(\frac{X_i}{\overline{\sigma}_i} - t\right)_+\bigg|\mathcal{F}_{i-1}\right] ~\le~ g(t),\quad\mbox{for all}\quad t \ge 0,
    \end{equation}
    where $\overline{\sigma}_i$ is an $\mathcal{F}_{i-1}$-measurable random variable. Additionally, set $\sigma_i^2 := \mathbb{E}[X_i^2|\mathcal{F}_{i-1}]$. Then for any $\delta\in(0, 1)$, $y\ge0$, and $v > 0$, with $Z\sim N(0, 1)$,
    \begin{equation}\label{eq:better-truncation-bound}
    \begin{split}
    \mathbb{P}\left(\bigcup_{k=1}^n\left\{\sum_{i=1}^k X_i \ge vQ_{5}\left(Z; \frac{\delta}{2}\right) + \frac{4v^2g(y)}{y\delta}\quad\mbox{and}\quad\sum_{i=1}^k \left(\frac{\overline{\sigma}_iy}{2} + \frac{\sigma_i^2}{2\overline{\sigma}_iy}\right) \le v^2\right\}\right)\le \delta.
    \end{split}
    \end{equation}
    Moreover, $Q_5(vZ; \delta/2) \le v\min\left\{\sqrt{2\log(2/\delta)},\, \Phi^{-1}(1 - \delta/11.4)\right\}.$
\end{corollary}
A proof of Corollary~\ref{cor:different-truncations} can be found in Section~\ref{appsec:different-truncations}.
Assumption~\eqref{eq:survival-assumption} is weaker than a tail assumption of the type $\mathbb{P}(X_i/\overline{\sigma}_i > t|\mathcal{F}_{i-1}) \le \overline{F}(t)$ for some survival function. Also, note that assumption~\eqref{eq:survival-assumption} is equivalent to
\begin{equation}\label{eq:positive-part-moment}
\mathbb{E}\left[\left((X_i/\overline{\sigma}_i)_+ - t\right)_+|\mathcal{F}_{i-1}\right] \le g(t)\quad\mbox{for all}\quad t \ge 0.
\end{equation}
See, for example, Theorem 2.1 of~\cite{marsiglietti2022moments} for cases where~\eqref{eq:positive-part-moment} might hold without the tail being dominated. 
Corollary~\ref{cor:different-truncations} can be thought of as an analogue of Fuk-Nagaev type inequalities~\citep{nagaev1979large}. Also, see~\cite{rio2017constants}. If $\mathbb{P}(X_i/\overline{\sigma}_i > u|\mathcal{F}_{i-1}) \le u^{-q}$ for some $q \ge 2$ and for all $u \ge 0$, then 
\[
\mathbb{E}[(X_i/\overline{\sigma}_i - t)_+|\mathcal{F}_{i-1}] \le \int_t^{\infty} \mathbb{P}(X_i/\overline{\sigma}_i > u|\mathcal{F}_{i-1})du = \int_t^{\infty} u^{-q}du = \frac{t^{-q + 1}}{q - 1}.
\]
In this case, inequality~\eqref{eq:better-truncation-bound} becomes
\begin{equation}\label{eq:better-tail-bound-finite-moment}
    \mathbb{P}\left(\bigcup_{k=1}^n \left\{\sum_{i=1}^k X_i \ge Q_5(vZ; \delta/2) + \frac{4v^2}{(q-1)\delta y^{q}}\quad\mbox{and}\quad \sum_{i=1}^k \left(\frac{\overline{\sigma}_iy}{2} + \frac{\sigma_i^2}{2\overline{\sigma}_iy}\right) \le v^2\right\}\right) \le \delta.
\end{equation}
This result resembles the conclusion of Corollary 4.2 of~\cite{rio2017constants} but with the leading term replaced by a Gaussian quantile rather than the crude upper bound typically obtained by the Cramer-Chernoff technique. However, it should be stressed that in our inequality $y$ is left as a free parameter. Unfortunately, this makes the bound weaker because changing $y$ with $\delta$ implies that the Gaussian leading term will no longer be a Gaussian term. 

In both inequalities~\eqref{eq:truncation-and-Thm2.1} and~\eqref{eq:better-truncation-bound}, improvements exist in the case where $\sigma_i^2$ are non-stochastic or bounded almost surely. For example, the results of~\cite{bentkus2010bounds} provide extensions of Theorem~\ref{thm:Freedman-improvement} with almost sure boundedness replaced by stochastic boundedness but requiring non-stochastic bound on the conditional variances. We could not combine the proof techniques in this note with the results of~\cite{bentkus2010bounds}, because we need a stochastic dominance for the ``worst-case'' random variables that appear in these bounds. In detail, in the context of improved Hoeffding and Bernstein inequalities, the worst-case random variables are two-point random variables~\eqref{eq:Hoeffding-bounded} and we know that sums of independent non-identically distributed Bernoulli's are stochastically dominated by sums of independent identically distributed Bernoulli with parameter being some aggregate of individual probabilities. The same holds with Gaussian and Poisson random variables, which lead to Theorem~\ref{thm:Azuma-Hoeffding}. 

Ideally, when using winsorization, one should control $\sum_{i=1}^n \min\{X_i, y\}$ or $\sum_{i=1}^n \min\{X_i, \overline{\sigma}_iy\}$ by using the fact that in addition to the conditional variance bounded by $\sigma_i^2$, these random variables also satisfy a tail assumption (governed by~\eqref{eq:survival-assumption}) before hitting $y$ as their upper bound. As stated currently, the inequalities~\eqref{eq:truncation-and-Thm2.1} and~\eqref{eq:better-truncation-bound} do not make use of this information, as was done in~\cite{rio2017constants} (for example). Using such tail information along with variance and bound, one could employ Theorem 3 of~\cite{bentkus2010bounds} for any martingale increment, but due to the lack of readily available stochastic domination of such sums, our bounds are weaker. For example, in the context of moment generating function (MGF) based proof in~\cite{rio2017constants}, one gets three regimes in the MGF for winsorized sums; see Proposition 3.5 and Eq. (3.14) there. These regimes show that the Gaussian leading term is unaffected by the winsorization (unlike ours), and the Poisson tail when balanced with winsorization yields the Fuk-Nagaev inequality with optimal constants. It might be possible to use Lemmas 4.1-4.8 of~\cite{pinelis2014bennett} or the results of~\cite{bentkus2008bounds} to obtain improvements of the Fuk-Nagaev inequalities, but our attempts to date have been unsuccessful. In summary, it would be interesting to explore Bentkus style bounds analogues to those in~\cite{Fan2015}.

In another direction, there are some results where truncation is introduced in the event of conditional variances. For example,~\cite{dzhaparidze2001bernstein} and~\cite{fan2017martingale} proved the following result for a supermartingale difference sequence $(X_i, \mathcal{F}_i)_{i\ge0}$:
\begin{align*}
&\mathbb{P}\left(\bigcup_{k=1}^n \left\{\sum_{i=1}^k X_i \ge x\quad\mbox{and}\quad \sum_{i=1}^k \mathbb{E}[X_i^2\mathbf{1}\{X_i \le y\}] + \sum_{i=1}^k X_i^2\mathbf{1}\{X_i > y\} \le v^2\right\}\right)\\ 
&\quad\le \inf_{\lambda \ge 0}\,\mathbb{E}\left[\exp(-\lambda x + \lambda y\widetilde{\Pi}_{v^2/y})\right].
\end{align*}
Although Theorem 2.1 of~\cite{fan2017martingale} is not stated in this form, it is equivalent; see Eq. (33) there. It stands to reason that a stronger inequality of the Bentkus type might also hold. In particular, we {\em conjecture} that the following holds
\begin{align*}
&\mathbb{P}\left(\bigcup_{k=1}^n \left\{\sum_{i=1}^k X_i \ge x\quad\mbox{and}\quad \sum_{i=1}^k \mathbb{E}[X_i^2\mathbf{1}\{X_i \le y\}] + \sum_{i=1}^k X_i^2\mathbf{1}\{X_i > y\} \le v^2\right\}\right)\\ 
&\quad\le \inf_{\lambda \ge 0}\,\mathbb{E}\left[\left(1 + \frac{\lambda}{2}(y\widetilde{\Pi}_{v^2/y} - x)\right)_+^2\right].
\end{align*}
It is not clear how to generalize the proof in~\cite{fan2017martingale} to the positive part polynomials.
\bibliography{notes}

\begin{thebibliography}{}

\bibitem[Bennett, 1962]{bennett1962probability}
Bennett, G. (1962).
\newblock Probability inequalities for the sum of independent random variables.
\newblock {\em Journal of the American Statistical Association},
  57(297):33--45.

\bibitem[Bentkus, 2002]{bentkus2002remark}
Bentkus, V. (2002).
\newblock A remark on {B}ernstein, {P}rokhorov, {B}ennett, {H}oeffding, and
  {T}alagrand inequalities.
\newblock {\em Lithuanian Mathematical Journal}, 42(3):262--269.

\bibitem[Bentkus, 2004]{bentkus2004hoeffding}
Bentkus, V. (2004).
\newblock On {H}oeffding’s inequalities.
\newblock {\em The Annals of Probability}, 32(2):1650--1673.

\bibitem[Bentkus, 2008]{bentkus2008extension}
Bentkus, V. (2008).
\newblock An extension of the {H}oeffding inequality to unbounded random
  variables.
\newblock {\em Lithuanian Mathematical Journal}, 48(2):137--157.

\bibitem[Bentkus, 2010]{bentkus2010bounds}
Bentkus, V. (2010).
\newblock Bounds for the stop loss premium for unbounded risks under the
  variance constraints.
\newblock {\em Preprint. Available at https://www.math.uni-bielefeld.
  de/sfb701/preprints/view/423}.

\bibitem[Bentkus and Ju{\v{s}}kevi{\v{c}}ius, 2008]{bentkus2008bounds}
Bentkus, V. and Ju{\v{s}}kevi{\v{c}}ius, T. (2008).
\newblock Bounds for tail probabilities of martingales using skewness and
  kurtosis.
\newblock {\em Lithuanian Mathematical Journal}, 48(1):30--37.

\bibitem[Dzhaparidze and Van~Zanten, 2001]{dzhaparidze2001bernstein}
Dzhaparidze, K. and Van~Zanten, J. (2001).
\newblock On bernstein-type inequalities for martingales.
\newblock {\em Stochastic processes and their applications}, 93(1):109--117.

\bibitem[Fan et~al., 2012]{fan2012hoeffding}
Fan, X., Grama, I., and Liu, Q. (2012).
\newblock Hoeffding’s inequality for supermartingales.
\newblock {\em Stochastic Processes and their Applications},
  122(10):3545--3559.

\bibitem[Fan et~al., 2015]{Fan2015}
Fan, X., Grama, I., and Liu, Q. (2015).
\newblock Exponential inequalities for martingales with applications.
\newblock {\em Electron. J. Probab.}, 20(1):1--22.

\bibitem[Fan et~al., 2017]{fan2017martingale}
Fan, X., Grama, I., and Liu, Q. (2017).
\newblock Martingale inequalities of type dzhaparidze and van zanten.
\newblock {\em Statistics}, 51(6):1200--1213.

\bibitem[Freedman, 1975]{freedman1975tail}
Freedman, D.~A. (1975).
\newblock On tail probabilities for martingales.
\newblock {\em the Annals of Probability}, pages 100--118.

\bibitem[Hahn and Klass, 1997]{hahn1997approximation}
Hahn, M.~G. and Klass, M.~J. (1997).
\newblock Approximation of partial sums of arbitrary iid random variables and
  the precision of the usual exponential upper bound.
\newblock {\em The Annals of Probability}, 25(3):1451--1470.

\bibitem[Hoeffding, 1963]{hoeffding1963probability}
Hoeffding, W. (1963).
\newblock Probability inequalities for sums of bounded random variables.
\newblock {\em Journal of the American Statistical Association},
  58(301):13--30.

\bibitem[Howard et~al., 2021]{howard2021time}
Howard, S.~R., Ramdas, A., McAuliffe, J., and Sekhon, J. (2021).
\newblock Time-uniform, nonparametric, nonasymptotic confidence sequences.
\newblock {\em The Annals of Statistics}, 49(2).

\bibitem[Kontorovich, 2014]{kontorovich2014concentration}
Kontorovich, A. (2014).
\newblock Concentration in unbounded metric spaces and algorithmic stability.
\newblock In {\em International conference on machine learning}, pages 28--36.
  PMLR.

\bibitem[Kutin, 2002]{kutin2002extensions}
Kutin, S. (2002).
\newblock Extensions to mcdiarmid’s inequality when differences are bounded
  with high probability.
\newblock {\em Dept. Comput. Sci., Univ. Chicago, Chicago, IL, USA, Tech. Rep.
  TR-2002-04}.

\bibitem[Lalley, 2014]{lalley2014martingale}
Lalley, S.~P. (2014).
\newblock The martingale central limit theorem.
\newblock Lecture notes, University of Chicago Department of Statistics.

\bibitem[Li and Liu, 2023]{li2023distribution}
Li, S. and Liu, Y. (2023).
\newblock Distribution-dependent mcdiarmid-type inequalities for functions of
  unbounded interaction.
\newblock In {\em International Conference on Machine Learning}, pages
  19789--19810. PMLR.

\bibitem[Li and Liu, 2024]{li2024concentration}
Li, S. and Liu, Y. (2024).
\newblock Concentration and moment inequalities for general functions of
  independent random variables with heavy tails.
\newblock {\em Journal of Machine Learning Research}, 25(268):1--33.

\bibitem[Marchina, 2018]{marchina2018concentration}
Marchina, A. (2018).
\newblock Concentration inequalities for separately convex functions.
\newblock {\em Bernoulli}, 24(4A):2906--2933.

\bibitem[Marsiglietti and Melbourne, 2022]{marsiglietti2022moments}
Marsiglietti, A. and Melbourne, J. (2022).
\newblock Moments, concentration, and entropy of log-concave distributions.
\newblock {\em arXiv preprint arXiv:2205.08293}.

\bibitem[Maurer, 2019]{maurer2019bernstein}
Maurer, A. (2019).
\newblock A bernstein-type inequality for functions of bounded interaction.
\newblock {\em Bernoulli}, 25(2):1451--1471.

\bibitem[Maurer and Pontil, 2021]{maurer2021concentration}
Maurer, A. and Pontil, M. (2021).
\newblock Concentration inequalities under sub-gaussian and sub-exponential
  conditions.
\newblock {\em Advances in Neural Information Processing Systems},
  34:7588--7597.

\bibitem[Nagaev, 1979]{nagaev1979large}
Nagaev, S.~V. (1979).
\newblock Large deviations of sums of independent random variables.
\newblock {\em The Annals of Probability}, pages 745--789.

\bibitem[Pinelis, 1998]{pinelis1998optimal}
Pinelis, I. (1998).
\newblock Optimal tail comparison based on comparison of moments.
\newblock In {\em High dimensional probability}, pages 297--314. Springer.

\bibitem[Pinelis, 2006a]{pinelis2006binomial}
Pinelis, I. (2006a).
\newblock Binomial upper bounds on generalized moments and tail probabilities
  of (super) martingales with differences bounded from above.
\newblock In {\em High dimensional probability}, pages 33--52. Institute of
  Mathematical Statistics.

\bibitem[Pinelis, 2006b]{Pinelis2006}
Pinelis, I. (2006b).
\newblock On normal domination of (super)martingales.
\newblock {\em Electron. J. Probab.}, 11:no. 39, 1049--1070.

\bibitem[Pinelis, 2014a]{pinelis2014bennett}
Pinelis, I. (2014a).
\newblock {On the Bennett-Hoeffding inequality}.
\newblock In {\em Annales de l'IHP Probabilit{\'e}s et statistiques},
  volume~50, pages 15--27.

\bibitem[Pinelis, 2014b]{pinelis2014optimal}
Pinelis, I. (2014b).
\newblock An optimal three-way stable and monotonic spectrum of bounds on
  quantiles: A spectrum of coherent measures of financial risk and economic
  inequality.
\newblock {\em Risks}, 2(3):349--392.

\bibitem[Rio, 2017]{rio2017constants}
Rio, E. (2017).
\newblock About the constants in the {F}uk-{N}agaev inequalities.
\newblock {\em Electronic Communications in Probability}, 22.

\bibitem[Talagrand, 1995]{talagrand1995missing}
Talagrand, M. (1995).
\newblock The missing factor in hoeffding's inequalities.
\newblock In {\em Annales de l'IHP Probabilit{\'e}s et statistiques},
  volume~31, pages 689--702.

\bibitem[van~de Geer, 2002]{van2002hoeffding}
van~de Geer, S.~A. (2002).
\newblock On hoeffding’s inequality for dependent random variables.
\newblock In {\em Empirical process techniques for dependent data}, pages
  161--169. Springer.

\end{thebibliography}
\bibliographystyle{apalike}
\newpage
\setcounter{section}{0}
\setcounter{equation}{0}
\setcounter{figure}{0}
\renewcommand{\thesection}{S.\arabic{section}}
\renewcommand{\theequation}{E.\arabic{equation}}
\renewcommand{\thefigure}{A.\arabic{figure}}
\renewcommand{\theHsection}{S.\arabic{section}}
\renewcommand{\theHequation}{E.\arabic{equation}}
\renewcommand{\theHfigure}{A.\arabic{figure}}
  \begin{center}
  \Large {\bf Supplement to ``On the Missing Factor in Some Concentration Inequalities for Martingales''}
  \end{center}
       
\begin{abstract}
This supplement contains the proofs of all the main results in the paper and some supporting lemmas. 
\end{abstract}
\section{Proof of Theorem~\ref{thm:Freedman-improvement}}\label{appsec:Freedman-improvement}
The proof proceeds as follows. We first provide the bound for $\mathbb{P}(\mathcal{E}_n)$ with $\mathcal{E}_n$ defined based on any supermartingale difference sequence satisfying the assumptions. Then we prove the bound for $\mathbb{P}(\mathcal{E}_k\mbox{ for some }1\le k\le n)$ by creating a new supermartingale difference sequence $(\widetilde{X}_i, \mathcal{F}_i)_{i\ge0}$ for which $\cup_{k=1}^n \mathcal{E}_k = \widetilde{\mathcal{E}}_n$, with $\widetilde{\mathcal{E}}_n$ defined using $\widetilde{X}_i, 1\le i\le n$. The proof of the inequality with $n=\infty$ follows from Proposition 4.8 of~\cite{pinelis2014bennett}.

Clearly, 
\begin{align*}
    \mathbb{P}(\mathcal{E}_n) &\le \mathbb{E}\left[\left(1 + \frac{\lambda}{2}\left(\sum_{i=1}^n X_i - x\right)\right)_+^2\mathbf{1}\{\mathcal{E}_n\}\right]\\
    &\le \mathbb{E}\left[\mathbb{E}\left[\left(1 + \frac{\lambda}{2}\left(\sum_{i=1}^n X_i - x\right)\right)_+^2\mathbf{1}\left\{\sum_{i=1}^n \sigma_i^2 \le v^2\right\}\bigg|\mathcal{F}_{n-1}\right]\right]
\end{align*}
Observe that
\[
\mathbb{E}[X_n|\mathcal{F}_{n-1}] \le 0,\quad \mbox{Var}(X_n|\mathcal{F}_{n-1}) = v^2 - \sum_{i=1}^{n-1}\sigma_i^2,\quad X_n \le 1\; \mbox{almost surely}.
\]
This implies that
\[
\mathbb{E}\left[\left(1 + \frac{\lambda}{2}\left(\sum_{i=1}^n X_i - x\right)\right)_+^2\mathbf{1}\{\mathcal{E}_n\}\bigg|\mathcal{F}_{n-1}\right] \le \mathbb{E}\left[\left(1 + \frac{\lambda}{2}\sum_{i=1}^{n-1} X_i + \frac{\lambda}{2}(G_n - x)\right)_+^2\mathbf{1}\left\{\sum_{i=1}^n \sigma_i^2 \le v^2\right\}\bigg|\mathcal{F}_{n-1}\right],
\]
where $G_n\in\{-(v^2 - \sum_{i=1}^{n-1}\sigma_i^2), 1\}$ with mean zero and variance bounded by $v^2 - \sum_{i=1}^{n-1}\sigma_i^2.$ Because $\sum_{i=1}^{n-1}\sigma_i^2$ is $\mathcal{F}_{n-2}$ measurable, we can think of $G_n$ (conditional on $\mathcal{F}_{n-2}$) to be a centered and scaled Bernoulli random variable. Let that (uncentered and unscaled) Bernoulli be $B_n$, independent of $\mathcal{F}_n$. Then there exists $\mathcal{F}_{n-2}$ measurable constants $a_n, b_n$ such that $G_n = a_nB_n + b_n$ and $\mathbb{E}[G_n|\mathcal{F}_{n-2}] = 0$.

Therefore,
\[
\mathbb{P}(\mathcal{E}_n) \le \mathbb{E}\left[\left(1 + \frac{\lambda}{2}\left(\sum_{i=1}^{n-1}X_i + G_n - x\right)\right)_+^2\right].
\]
Now consider the conditional expectation given $\mathcal{F}_{n-2}\cup\{B_n\}$. We know
\[
\mathbb{E}[X_{n-1}|\mathcal{F}_{n-2}\cup\{B_n\}] \le 0,\quad\mbox{Var}(X_{n-1}|\mathcal{F}_{n-2}\cup\{B_n\}) = \sigma_{n-1}^2,\quad X_{n-1} \le 1\quad\mbox{almost surely.}
\]
This implies
\[
\mathbb{E}[(1 + (\lambda/2)(X_{n-1} - u))_+^2|\mathcal{F}_{n-2}] \le \mathbb{E}[(1 + (\lambda/2)(G_{n-1} - u))_+^2|\mathcal{F}_{n-2}],\quad\mbox{for all}\quad \lambda \ge 0, u\in\mathbb{R}.
\]
Here $G_{n-1}\in\{-\sigma_{n-1}^2, 1\}$ with mean zero (conditional on $\mathcal{F}_{n-2}$) and variance $\sigma_{n-1}^2$. As before, we can write $G_{n-1} = a_{n-1}B_{n-1} + b_{n-1}$ for some $\mathcal{F}_{n-2}$ measurable random variables $a_{n-1}, b_{n-1}$ and Bernoulli random variable $B_{n-1}$ independent of $\mathcal{F}_n$ so that $\mathbb{E}[G_{n-1}|\mathcal{F}_{n-2}] = 0$. 
Therefore, 
\[
\mathbb{P}(\mathcal{E}_n) \le \mathbb{E}\left[\left(1 + \frac{\lambda}{2}\left(\sum_{i=1}^{n-2}X_i + G_{n-1} + G_n - x\right)\right)_+^2\right],
\]
Note that conditional on $\mathcal{F}_{n-2}$, $G_{n}$ and $G_{n-1}$ are two mean zero independent random variables and $\mbox{Var}(G_n + G_{n-1}|\mathcal{F}_{n-2}) = \sigma_{n-1}^2 + v^2 - \sum_{i=1}^{n-1}\sigma_{i}^2 = v^2 - \sum_{i=1}^{n-2} \sigma_i^2$. Hence, by Lemma 4.5 of~\cite{bentkus2004hoeffding}, we can bound the right hand side with
\[
\mathbb{P}(\mathcal{E}_n) \le \mathbb{E}\left[\left(1 + \frac{\lambda}{2}\left(\sum_{i=1}^{n-2}X_i + \tilde{G}_{n-1} + \tilde{G}_n - x\right)\right)_+^2\right],
\]
where $\tilde{G}_{n-1}, \tilde{G}_n$ are IID conditional on $\mathcal{F}_{n-2}$ and have mean zero, bounded above by 1, and equal variances of $(v^2 - \sum_{i=1}^{n-2}\sigma_i^2)/2$. Note that this variance is $\mathcal{F}_{n-3}$-measurable, which implies that we can take $\tilde{G}_{n-1}, \tilde{G}_n$ to be IID conditional on $\mathcal{F}_{n-3}$.

Now consider this expectation conditional on $\mathcal{F}_{n-3}\cup\{\tilde{B}_{n-1}, \tilde{B}_n\}$. We have
\[
\mathbb{E}[X_{n-2}|\mathcal{F}_{n-3}] \le 0, \mbox{Var}(X_{n-2}|\mathcal{F}_{n-3}) \le \sigma_{n-2}^2,\mbox{ and }X_{n-3} \le 1\mbox{ almost surely.}
\]
Hence, we can replace $X_{n-2}$ with $G_{n-2}$. Observe that conditional on $\mathcal{F}_{n-3}$, $G_{n-2}, \tilde{G}_{n-1}, \tilde{G}_{n}$ are independent random variables taking only two values, and $\mbox{Var}(G_{n-2} + \tilde{G}_{n-1} + \tilde{G}_{n}|\mathcal{F}_{n-3}) = v^2 - \sum_{i=1}^{n-3}\sigma_i^2$. Applying Lemma 4.5 of Bentkus (2004) again and iterating this procedure, we can bound the probability by
\[
\mathbb{P}(\mathcal{E}_n) \le \mathbb{E}\left[\left(1 + \frac{\lambda}{2}\left(\sum_{i=1}^n \tilde{G}_i - x\right)\right)_+^2\right],
\]
where $\tilde{G}_i, 1\le i\le n$ are IID binary random variables taking values in $\{-v^2/n, 1\}$ with mean zero and variance $v^2/n$.

To bound the probability of $\cup_{k=1}^n \mathcal{E}_k$, define
$T\in\{1, 2, \ldots, n\}$ be the minimum $k$ such that $\sum_{i=1}^k X_i \ge x$ and $\sum_{i=1}^k \sigma_i^2 \le v^2$. If $\cup_{k=1}^n \mathcal{E}_k$ does not occur, set $T = n$. This is a stopping time because $\{T \le i\}$ can be determined by looking at $\mathcal{F}_i$, or equivalently, $\{T \le i\}\in\mathcal{F}_i$. Then 
\[
\mathbf{1}\{\mathcal{E}_n\} \le \mathbf{1}\left\{\sum_{i=1}^n \widetilde{X}_i \ge x,\quad\mbox{and}\quad \sum_{i=1}^n \widetilde{\sigma}_i^2 \le v^2\right\},
\]
where $\widetilde{X}_i = X_i\mathbf{1}\{T \ge i\}$ and $\widetilde{\sigma}_i^2 = \mathbb{E}[\widetilde{X}_i^2|\mathcal{F}_{i-1}]$. Observe that $\sum_{i=1}^n \widetilde{X}_i = \sum_{i=1}^T X_i$. Also, note that $(\widetilde{X}_k,\mathcal{F}_k)_{k\ge0}$ is a supermartingale difference sequence if $(X_k, \mathcal{F}_k)$ is a supermartingale difference sequence. This follows because 
\[
\mathbb{E}\left[\sum_{i=1}^{k+1} \tilde{X}_i\big|\mathcal{F}_{k}\right] = \sum_{i=1}^{k+1}\mathbb{E}[X_i\mathbf{1}\{T \ge i\}\big|\mathcal{F}_k] = \sum_{i=1}^k X_i\mathbf{1}\{T \ge i\} + \mathbb{E}[X_{k+1}\mathbf{1}\{T \ge k+1\}\big|\mathcal{F}_k],
\]
where the last equality follows because $X_i\in\mathcal{F}_k$ for any $k \ge i$ and $\mathbf{1}\{T \ge i\} = 1 - \mathbf{1}\{T \le i-1\} \in \mathcal{F}_{i-1} \subseteq \mathcal{F}_{k}$ for any $k \ge i$. Finally, because $\mathbf{1}\{T \ge k+1\}\in\mathcal{F}_k$, we get $\mathbb{E}[X_{k+1}\mathbf{1}\{T \ge k\}|\mathcal{F}_k] = \mathbf{1}\{T \ge k\}\mathbb{E}[X_{k+1}|\mathcal{F}_k] \le 0$. This shows that if $(X_k, \mathcal{F}_k)_{k\ge 0}$ is a supermartingale, then $(\widetilde{X}_k, \mathcal{F}_k)_{k\ge 0}$ is also a supermartingale. Hence, applying the previous result, we get
\[
\mathbb{P}\left(\bigcup_{k=1}^n \mathcal{E}_k\right) \le \inf_{\lambda \ge 0}\,\mathbb{E}\left[\left(1 + \frac{\lambda}{2}\left(\sum_{i=1}^n {G}_i - x\right)\right)_+^2\right],
\]
where ${G}_i, 1\le i\le n$ are IID mean zero random variables such that $\mathbb{P}({G_i} \in \{-v^2/n, 1\}) = 1$.

Finally, from Proposition 2.8 of~\cite{pinelis2014bennett}, it follows that for any $\lambda \ge 0$,
\[
\mathbb{E}\left[\left(1 + \frac{\lambda}{2}\left(\sum_{i=1}^n \widetilde{G}_i - x\right)\right)_+^2\right] \le \mathbb{E}\left[\left(1 + \frac{\lambda}{2}(\widetilde{\Pi}_{v^2} - x)\right)_+^2\right],
\]
where $\widetilde{\Pi}_{v^2}$ is a centered Poisson$(v^2)$ random variable. Given that the right hand side is independent of $n$, we get
\[
\mathbb{P}\left(\sum_{i=1}^k X_i \ge x\;\mbox{and}\; \sum_{i=1}^k \sigma_i^2 \le v^2\;\mbox{for some}\; k\ge 1\right) \le \inf_{\lambda \ge 0}\, \mathbb{E}\left[\left(1 + \frac{\lambda}{2}(\widetilde{\Pi}_{v^2} - x)\right)_+^2\right] \le \frac{e^2}{2}\mathbb{P}^\circ(\widetilde{\Pi}_{v^2} \ge x),
\]
where $x\mapsto\mathbb{P}^\circ(\widetilde{\Pi}_{v^2} \ge x)$ is the log-concave majorant of $x\mapsto \mathbb{P}(\widetilde{\Pi}_{v^2} \ge x)$. 
\section{Proof of Theorem~\ref{thm:Azuma-Hoeffding}}\label{appsec:proof-of-Azuma-Hoeffding}
We have, as before,
\[
\mathbb{P}(\mathcal{E}_n) \le \mathbb{E}\left[\mathbb{E}\left[\left(1 + \frac{\lambda}{5}\left(\sum_{i=1}^n X_i - x\right)\right)_+^5\bigg|\mathcal{F}_{n-1}\right]\mathbf{1}\left\{s_n^2 \le v^2 - V_{n-1}^2\right\}\right].
\]
Conditional on $\mathcal{F}_{n-1}$, $X_n$ has a non-positive conditional mean, bounded by $B_n$, and conditional variance bounded by $\sigma_n^2$ so that $(1/2)(B_n + \sigma_n^2/B_n) \le (v^2 - V_{n-1}^2)^{1/2}.$ Applying Lemma 5.1.3 of~\cite{Pinelis2006}, we conclude
\begin{align*}
&\mathbb{E}\left[\left(1 + \frac{\lambda}{5}\left(\sum_{i=1}^n X_i - x\right)\right)_+^5\bigg|\mathcal{F}_{n-1}\right]\mathbf{1}\left\{s_n^2 \le v^2 - V_{n-1}^2\right\}\\ 
&\quad\le \mathbb{E}\left[\left(1 + \frac{\lambda}{5}\left(\sum_{i=1}^{n-1} X_i + G_n - x\right)\right)_+^5\bigg|\mathcal{F}_{n-1}\right]\mathbf{1}\{V_{n-1} \le v^2\},
\end{align*}
where $G_n|\mathcal{F}_{n-1}\sim N(0, v^2 - V_{n-1}^2)$. Observe that $v^2 - V_{n-1}^2$ is $\mathcal{F}_{n-2}$-measurable, which implies that $G_n|\mathcal{F}_{n-2} \sim N(0, v^2 - V_{n-1}^2)$. This implies
\begin{align*}
\mathbb{P}(\mathcal{E}_n) &\le \mathbb{E}\left[\left(1 + \frac{\lambda}{5}\left(\sum_{i=1}^{n-1}X_i + G_n - x\right)\right)_+^5\mathbf{1}\{V_{n-1}^2 \le v^2\}\right]\\
&\le \mathbb{E}\left[\mathbb{E}\left[\left(1 + \frac{\lambda}{5}\left(\sum_{i=1}^{n-2}X_i + G_n + X_{n-1} - x\right)\right)_+^5\mathbf{1}\{V_{n-1}^2 \le v^2\}\bigg|\mathcal{F}_{n-2}, G_n\right]\right].
\end{align*}
Observe now that conditional on $\mathcal{F}_{n-2}, G_n$, $X_{n-1}$ has a non-positive mean, variance $\sigma_{n-1}^2$, bounded by $B_{n-1}$ so that $(1/2)(B_{n-1} + \sigma_{n-1}^2/B_{n-1}) \le s_{n-1} \in \mathcal{F}_{n-2}$. Therefore, by Lemma 5.1.3 of~\cite{Pinelis2006}, we conclude that on the event $V_{n-1}^2 \le v^2$,
\begin{align*}
&\mathbb{E}\left[\left(1 + \frac{\lambda}{5}\left(\sum_{i=1}^{n-2}X_i + G_n + X_{n-1} - x\right)\right)_+^5\bigg|\mathcal{F}_{n-2}, G_n\right]\\ 
&\quad\le \mathbb{E}\left[\left(1 + \frac{\lambda}{5}\left(\sum_{i=1}^{n-2} X_i + G_n + G_{n-1}\right)\right)_+^5\bigg|\mathcal{F}_{n-2}, G_n\right],
\end{align*}
where $G_{n-1}|\mathcal{F}_{n-2}, G_n \sim N(0, s_{n-1}^2)$. This implies that $G_{n-1} + G_n|\mathcal{F}_{n-2} \sim N(0, v^2 - V_{n-1}^2 + s_{n-1}^2) = N(0, v^2 - V_{n-2}^2)$. Combining, we get
\[
\mathbb{P}(\mathcal{E}_n) \le \mathbb{E}\left[\left(1 + \frac{\lambda}{5}\left(\sum_{i=1}^{n-2} X_i + (v^2 - V_{n-2}^2)^{1/2}Z - x\right)\right)_+^5\mathbf{1}\{V_{n-2}^2 \le v^2\}\right],
\]
where $Z\sim N(0, 1)$ independent of $\mathcal{F}_n$.
Repeating the argument iteratively, we conclude that
\[
\mathbb{P}(\mathcal{E}_n) \le \mathbb{E}\left[\left(1 + \frac{\lambda}{5}(vZ - x)\right)_+^5\right].
\]
Hence, taking the infimum over $\lambda \ge 0$, we get
\[
\mathbb{P}(\mathcal{E}_n) \le \inf_{\lambda \ge 0}\mathbb{E}\left[\left(1 + \frac{\lambda}{5}(vZ - x)\right)_+^5\right] \le  5!(e/5)^5\mathbb{P}\left(vZ \ge x\right).
\]

Now consider the bigger event
\[
\widetilde{\mathcal{E}}_n = \left\{\sum_{i=1}^k X_i \ge x\quad\mbox{and}\quad V_k^2 \le v^2\mbox{ for some }1 \le k\le n\right\}.
\]
As before, define $T$ as the smallest $k\in\{1, 2, \ldots, n\}$ such that $\sum_{i=1}^k X_i \ge x$ and $V_k^2 \le v^2$. If no such $k$ exists, then set $T = n$. This is a stopping time, and set $\widetilde{X}_i = X_i\mathbf{1}\{T \ge i\}$. Then it is as before easy to verify that $(\widetilde{X}_i, \mathcal{F}_i), 1 \le i\le n$ is a supermartingale if $(X_i, \mathcal{F}_i)$ is a supermartingale. Moreover, $\widetilde{\sigma}_i^2 = \mbox{Var}(\widetilde{X}_i|\mathcal{F}_{i-1}) = \sigma_i^2\mathbf{1}\{T \ge i\}$ and set $\widetilde{s}_i = (1/2)(B_i + \widetilde{\sigma}_i^2/B_i)$. Hence, 
\[
\mathbf{1}\{\widetilde{\mathcal{E}}_n\} \le \mathbf{1}\left\{\sum_{i=1}^n \widetilde{X}_i \ge x\mbox{ and }\sum_{i=1}^n \widetilde{s}_i^2 \le v^2\right\}.
\]
and
\[
\mathbb{P}(\widetilde{\mathcal{E}}_n) \le \inf_{\lambda \ge 0}\,\mathbb{E}\left[\left(1 
+ \frac{\lambda}{5}(vZ - x)\right)_+^5\right].
\]
This completes the proof.
\section{Proof of Lemma~\ref{lem:Hoeffding-unbounded}}\label{appsec:proof-of-Hoeffding-unbounded}
\begin{lemma}\label{lem:inequalities_from_dominance}
    Suppose $X$ and $Y$ are two real-valued random variables satisfying 
    \begin{equation}\label{eq:order_1_stochastic}
    \mathbb{E}[(X - t)_+] \le \mathbb{E}[(Y - t)_+]\quad\mbox{for all}\quad t\in\mathbb{R}.
    \end{equation}
    Then $\mathbb{E}[X] \le \mathbb{E}[Y]$. In other words, $X\preceq_1 Y$ implies $\mathbb{E}[X] \le \mathbb{E}[Y]$.
\end{lemma}
\begin{proof}
The fact that $X\preceq_1 Y$ implies $\mathbb{E}[X] \le \mathbb{E}[Y]$ follows from Proposition 4 of Bentkus (2008, Lithuanian Mathematical Journal, 48(3)). The following facts will be important: $x = x_+ - x_ = x_+ - (-x)_+$.
The proof is short and is as follows: Observe that
\[
\mathbb{E}[X - t] = \mathbb{E}[(X - t)_+] - \mathbb{E}[(t - X)_+] \le \mathbb{E}[(Y - t)_+] = \mathbb{E}[Y - t] + \mathbb{E}[(t - Y)_+].
\]
Therefore, 
\[
\mathbb{E}[X] \le \mathbb{E}[Y] + \mathbb{E}[(t - Y)_+]\quad\mbox{for all}\quad t\in\mathbb{R}.
\]
Taking the limiting as $t\to-\infty$ implies the result.
\end{proof}
\begin{proof}[Proof of Lemma~\ref{lem:Hoeffding-unbounded}]
    The inequality $\mathbb{E}[T] \le \mathbb{E}[X] \le \mathbb{E}[W]$ follows from Lemma~\ref{lem:inequalities_from_dominance}. Because $T \preceq_0 W$ (i.e., $\mathbb{P}(T > x) \le \mathbb{P}(W > x)$ for all $x$), we get $a_q \le b_q$ for all $q\in[0, 1]$. If $q = 0$, then $a_q$ is the right end of the support of $T$ and hence, $\mathbb{P}(\xi_q \le x) = \mathbb{P}(T \le x)$ for all $x < a_q$ and $\mathbb{P}(\xi_q \le a_q) = 1$. This implies $\xi_0 \overset{d}{=} T$. Similarly, $\xi_1 \overset{d}{=} W$. 

    Setting $F_{\xi_q}(x) = \mathbb{P}(\xi_q \le x)$, $F_T(x) = \mathbb{P}(T \le x),$ and $F_W(x) = \mathbb{P}(W \le x)$, we get
    \[
    \mathbb{E}[\xi_q] = \int_0^1 F_{\xi_q}^{-1}(\delta)d\delta = \int_0^{1-q} F_T^{-1}(\delta)d\delta + \int_{1-q}^1 F_W^{-1}(\delta)d\delta = \mathbb{E}[\eta] + \int_{1-q}^1 (F_W^{-1}(\delta) - F_T^{-1}(\delta))d\delta. 
    \]
    Because $T \preceq_0 W$, we get differentiability and non-decreasing property of $q\mapsto \mathbb{E}[\xi_q]$.

    To show the final property $X\preceq_1 \xi_{q_0}$, we consider three cases: (i) $t \le a_{q_0}$, (ii) $a_{q_0} < t < b_{q_0}$, (iii) $t \ge b_{q_0}$. The first and third cases are easy to prove from $-X \preceq_1 -T$ and $X \preceq_1 W$. For example, if $t \le a_{q_0}$, 
    \[
    \mathbb{E}[(X - t)_+] = \mathbb{E}[(X - t)] + \mathbb{E}[(t - X)_+] \le \mathbb{E}[(\xi_{q_0} - t)] + \mathbb{E}[(t - T)_+].
    \]
    For $t \le a_{q_0}$,
    \[
    \mathbb{E}[(t - T)_+] = \int_0^{\infty} \mathbb{P}(t - T > s)ds = \int_{-\infty}^t \mathbb{P}(T < u)du = \int_{-\infty}^t \mathbb{P}(\xi_{q_0} < u)du = \mathbb{E}[(t - \xi_{q_0})_+],
    \]
    where the condition $t \le a_{q_0}$ is used in the penultimate equality. Therefore, $\mathbb{E}[(X - t)_+] \le \mathbb{E}[(\xi_{q_0} - t)_+]$ for $t \le a_{q_0}$.

    If $a_{q_0} < t < b_{q_0}$, then note that
    \[
    \mathbb{E}[(\xi_{q_0} - t)_+] = \int_{t}^{b_{q_0}} \mathbb{P}(\xi_{q_0} > s)ds + \int_{b_{q_0}}^{\infty} \mathbb{P}(\xi_{q_0} > s)ds = q_0(b_{q_0} - t) + \int_{b_{q_0}}^{\infty} \mathbb{P}(W > s)ds = q_0(b_{q_0} - t) + \mathbb{E}[(W - t)_+].
    \]
    On the other hand, we can write
    \[
    t = \frac{b_{q_0} - t}{b_{q_0} - a_{q_0}}a_{q_0} + \frac{t - a_{q_0}}{b_{q_0} - a_{q_0}}b_{q_0},
    \]
    which implies
    \[
    \mathbb{E}[(X - t)_+] \le \frac{b_{q_0} - t}{b_{q_0} - a_{q_0}}\mathbb{E}[(X - a_{q_0})_+] + \frac{t - a_{q_0}}{b_{q_0} - a_{q_0}}\mathbb{E}[(X - b_{q_0})_+].
    \]
    From case (i) \& (iii), this can be upper bounded as
    \[
    \mathbb{E}[(X - t)_+] \le \frac{b_{q_0} - t}{b_{q_0} - a_{q_0}}\mathbb{E}[(\xi_{q_0} - a_{q_0})_+] + \frac{t - a_{q_0}}{b_{q_{0}} - a_{q_0}}\mathbb{E}[(W - b_{q_0})_+] = \mathbb{E}[(\xi_{q_0} - t)_+],
    \]
    where the last equality follows from the expression for $\mathbb{E}[(\xi_{q_0} - t)_+]$.
\end{proof}
\section{Proof of Corollary~\ref{cor:different-truncations}}\label{appsec:different-truncations}
    For any $\delta\in(0, 1)$, set
    \begin{equation}\label{eq:quantile-thresholds}
    x_1 = \inf_{t\in\mathbb{R}}\, t + \frac{(\mathbb{E}[(vZ - t)_+])^{1/5}}{(\delta/2)^5},\quad\mbox{and}\quad x_2 = \frac{4v^2g(y)}{y\delta}.
    \end{equation}
    From the proof of Theorem~\ref{thm:Freedman-improvement}, it suffices to prove that $\mathbb{P}(\mathcal{E}_n) \le \delta$,
    where
    \[
    \mathcal{E}_n = \left\{\sum_{i=1}^n X_i \ge x_1 + x_2\quad\mbox{and}\quad\sum_{i=1}^n \left(\frac{\overline{\sigma}_iy}{2} + \frac{\sigma_i^2}{2\overline{\sigma}_iy}\right) \le v^2\right\}.
    \]
    Observe that $\mathcal{E}_n \subseteq \mathcal{E}_n^{(1)} \cup \mathcal{E}_n^{(2)}$ where
    \begin{align*}
        \mathcal{E}_n^{(1)} &:= \left\{\sum_{i=1}^n \min\{X_i, \overline{\sigma}_iy\} \ge x_1\;\mbox{and}\;\sum_{i=1}^n \left(\frac{\overline{\sigma}_iy}{2} + \frac{\sigma_i^2}{2\overline{\sigma}_iy}\right) \le v^2\right\},\\
        \mathcal{E}_n^{(2)} &:= \left\{\sum_{i=1}^n (X_i - \overline{\sigma}_iy)_+ \ge x_2\;\mbox{and}\; \sum_{i=1}^n \overline{\sigma}_i \le \frac{2v^2}{y}\right\}.
    \end{align*}
    Because $(\min\{X_i, \overline{\sigma}_iy\}, \mathcal{F}_i)_{i\ge0}$ is a supermartingale difference sequence, Theorem~\ref{thm:Azuma-Hoeffding} implies that 
    \begin{equation}\label{eq:event-1}
    \mathbb{P}(\mathcal{E}_n^{(1)}) \le \inf_{\lambda \ge 0}\,\mathbb{E}\left[\left(1 + \frac{\lambda}{5}(vZ - x_1)\right)_+^5\right] ~\le~ \min\left\{e^{-x_1^2/(2v^2)},\, \frac{5!e^5}{5^5}(1 - \Phi(x_1/v))\right\}.
    \end{equation}
    From Eq. (2.2), Theorem 2.3 of~\cite{pinelis2014optimal}, and our choice of $x_1$ in~\eqref{eq:quantile-thresholds}, we get $\mathbb{P}(\mathcal{E}_n^{(1)}) \le \delta/2$. The second inequality of~\eqref{eq:event-1} implies that 
    \[
    Q_5(vZ; \delta/2) \le \min\left\{v\sqrt{2\log(2/\delta)},\, v\Phi^{-1}(1 - \delta/11.4)\right\}.
    \]
    
    To bound the probability of $\mathcal{E}_n^{(2)}$, set $\mathcal{V}_k := \mathbf{1}\{\sum_{i=1}^k \overline{\sigma}_i \le 2v^2/y\}$ for $1\le k\le n$. Note that, by assumption, $\overline{\sigma}_i\in\mathcal{F}_{i-1}$, $\mathcal{V}_k\in\mathcal{F}_{k-1}$. This fact will be used in the following without specific reference.
    Note that, by Markov inequality, 
    \begin{align*}
    \mathbb{P}(\mathcal{E}_n^{(2)}) &\le \frac{1}{x_2}\mathbb{E}\left[\sum_{i=1}^n (X_i - \overline{\sigma}_iy)_+\mathbf{1}\left\{\mathcal{V}_n\right\}\right]\\
    &\le \frac{1}{x_2}\mathbb{E}\left[\sum_{i=1}^{n-1} (X_i - \overline{\sigma}_iy)_+\mathbf{1}\{\mathcal{V}_{n-1}\} + \left(\frac{2v^2}{y} - \sum_{i=1}^{n-1} \overline{\sigma}_i\right)\mathbf{1}\{\mathcal{V}_{n-1}\}\mathbb{E}[(X_n/\overline{\sigma}_n - y)_+|\mathcal{F}_{n-1}]\right]\\
    &\le \frac{1}{x_2}\mathbb{E}\left[\sum_{i=1}^{n-1} (X_i - \overline{\sigma_i}y)_+\mathbf{1}\{\mathcal{V}_{n-1}\} + \left(\frac{2v^2}{y} - \sum_{i=1}^{n-1} \overline{\sigma}_i\right)g(y)\mathbf{1}\{\mathcal{V}_{n-1}\}\right].
    \end{align*}
    This implies, again by~\eqref{eq:survival-assumption},
    \begin{align*}
    \mathbb{P}(\mathcal{E}_n^{(2)}) &\le \frac{1}{x_2}\mathbb{E}\left[\sum_{i=1}^{n-2}(X_i - \overline{\sigma}_iy)_+\mathbf{1}\{\mathcal{V}_{n-2}\} + \left(\frac{2v^2}{y} - \sum_{i=1}^{n-1} \overline{\sigma}_i\right)\mathbf{1}\{\mathcal{V}_{n-1}\}g(y)\right]\\ 
    &\quad+ \mathbb{E}\left[\overline{\sigma}_{n-1}\mathbf{1}\{\mathcal{V}_{n-1}\}\mathbb{E}[(X_{n-1}/\overline{\sigma}_{n-1} - y)_+|\mathcal{F}_{n-2}]\right]\\
    &\le \frac{1}{x_2}\mathbb{E}\left[\sum_{i=1}^{n-2} (X_i - \overline{\sigma}_iy)_+\mathbf{1}\{\mathcal{V}_{n-2}\} + \left(\frac{2v^2}{y} - \sum_{i=1}^{n-2} \overline{\sigma}_i\right)\mathbf{1}\{\mathcal{V}_{n-2}\}g(y)\right].
    \end{align*}
    Repeating this argument further, we get
    \[
    \mathbb{P}(\mathcal{E}_n^{(2)}) \le \frac{2v^2g(y)}{yx_2}.
    \]
    Our choice of $x_2$ in~\eqref{eq:quantile-thresholds} implies that $\mathbb{P}(\mathcal{E}_n^{(2)}) \le \delta/2$. This completes the proof of the inequality $\mathbb{P}(\mathcal{E}_n) \le \delta$ and proves the result.
\end{document}